\documentclass[11pt]{amsart}
\usepackage{amsmath,amssymb,mathrsfs,amscd}
\usepackage{graphicx}
\usepackage[all]{xy}

\thanks
{$^{\dag}$Partially supported by
Research Fellowships of the 
Japan Society for the Promotion
of Science for Young Scientists.
}

\title{Examples of compact minitwistor spaces and their moduli space}
\author{Nobuhiro Honda
}  

\newcommand{\ol}{\overline}
\newcommand{\ra}{\rightarrow}

\newcommand{\set}{\,|\,}
\newcommand{\proofend}{\hfill$\square$}

\newcommand{\vsp}{\vspace{3mm}}

\newtheorem{prop}{Proposition}[section]
\newtheorem{lemma}[prop]{Lemma}

\newtheorem{thm}[prop]{Theorem}

\begin{document}

\maketitle

\begin{abstract} 
In a paper \cite{H-bitan} we obtained explicit examples of Moishezon twistor spaces of some  compact self-dual four-manifolds admitting a non-trivial Killing field, and also determined their moduli space.
In this note we investigate minitwistor spaces associated to these twistor spaces.
We determine their structure, minitwistor lines and also their moduli space, by using a  double covering structure of the twistor spaces.
In particular, we find that these minitwistor spaces have different properties in many respects, compared to known examples of minitwistor spaces.
Especially, we show that the moduli space of the minitwistor spaces is identified with the configuration space of different 4 points on a circle divided by the standard ${\rm{PSL}}\,(2,\mathbf R)$-action.
\end{abstract}



\section{Introduction} 
In \cite{JT85} P.E. Jones and K.P. Tod established a reduction theory for  self-dual 4-manifolds with a non-trivial Killing field.
We briefly recall their results.
Suppose that a self-dual metric $g$ on a 4-manifold $M$ admits a free isometric ${\rm{U}}(1)$-action. Then the quotient 3-manifold $M/{\rm{U}}(1)$ is naturally equipped with so called a {\em Weyl structure}, which is 
a pair of a conformal structure (associated to the natural  Riemannian metric on $M/{\rm{U}}(1)$) and an affine connection compatible with the conformal structure.
As a consequence of the self-duality of  $g$, a curvature of the affine connection satisfies a kind of Einstein condition and the pair becomes
{\em Einstein-Weyl structure} in the sense of N.\,J.\,Hitchin \cite{Hi82-1}.
Moreover, the function on $M/{\rm{U}}(1)$ obtained by associating the length of each ${\rm{U}}(1)$-orbits (with respect to $g$) satisfies certain linear equation, which is called a {\em monopole equation}. 
Thus, the Einstein-Weyl condition and the monopole equation can be thought as a non-linear  and linear part of the self-duality equation respectively.
This construction is invertible. Namely, if a 3-manifold $N$ is equipped with an Einstein-Weyl structure, and if $M\ra N$ is a principal ${\rm{U}}(1)$-bundle equipped with a (positive) solution of the monopole equation, a conformal structure on $M$ is naturally constructed and it becomes self-dual. One can also refer \cite{LB93} for details.

The last inversion of the reduction theory already produces non-trivial self-dual metrics even if one takes the flat Euclidean 3-space (with the natural Einstein-Weyl structure) and if one allows a certain kind of singularities for the ${\rm{U}}(1)$-bundle. 
Namely, the Eguchi-Hanson metric \cite{EH79} on the cotangent bundle of $\mathbf{CP}^1$ and the Gibbons-Hawking metrics \cite{GH78} on the minimal resolution of $\mathbf C^2/\Gamma$, where $\Gamma$ is a cyclic subgroup of $\rm{SU}(2)$,  are obtained in this way.
Later, C. LeBrun \cite{LB91} successfully applies the inversion construction for the hyperbolic 3-space  (again with the natural Einstein-Weyl structure) to realize  explicit self-dual metrics on $n\mathbf{CP}^2$, the connected sum of $n$ copies of complex projective planes.

This reduction theory for self-dual metrics can be translated into that of twistor spaces \cite{JT85} (see also \cite{LB91,LB93}). In particular, taking the quotient of a self-dual manifold by  a ${\rm{U}}(1)$-action corresponds to taking the quotient of the associated twistor space by the natural holomorphic $\mathbf C^*$-action. 
Because  $\mathbf C^*$-action can become pathological in general as seen in \cite[\S 4]{Hon01},  the quotient space of the $\mathbf C^*$-action does not necessarily possess a structure of a complex surface.
But  if the twistor space is  Moishezon for example, such a pathology does not occur and the quotient space can have a natural structure of a complex surface (with singularities in general). 
In this case, the quotient complex surface is called a {\em minitwistor space}.
For the Gibbons-Hawking metrics and the LeBrun metrics mentioned above, this is indeed the case;  the resulting minitwistor spaces are the total space of the holomorphic tangent bundle of $\mathbf{CP}^1$ for the Gibbons-Hawking metrics, and the product $\mathbf{CP}^1\times\mathbf{CP}^1$ for the LeBrun metrics.
Here, an important feature in these two basic examples is that, while the self-dual (or hyper-K\"ahler) structure actually deforms, the corresponding Einstein-Weyl 3-manifolds, and hence their minitwistor spaces, do not deform.
 
In a paper \cite{H-bitan} the author explicitly constructed a family of twistor spaces on $3\mathbf{CP}^2$ parametrized by a 3-dimensional connected space.
The corresponding self-dual metrics have a non-trivial ${\rm{U}}(1)$-action but are not conformal to the LeBrun metrics.
Moreover, the constructed family is complete in the sense that, on $3\mathbf{CP}^2$, every non-LeBrun self-dual metric with ${\rm{U}}(1)$-action is a member of this family, at least if the self-dual metric is supposed to have a positive scalar curvature.
The purpose of this note is to investigate minitwistor spaces of these explicit twistor spaces of $3\mathbf{CP}^2$.
Our first result is a determination of the structure of these minitwistor spaces; namely we show that the minitwistor spaces have a natural structure of a branched double covering of $\ol{\Sigma}_2$, the Hirzebruch surface of degree 2 with the $(-2)$-section contracted, and that the branching divisor is a smooth anticanonical curve which is  an elliptic curve disjoint from the node of $\ol{\Sigma}_2$ (Theorem \ref{thm-main1}).
Briefly speaking this is a reflection of the property that our twistor spaces of $3\mathbf{CP}^2$ have a structure of generically 2 to 1 covering of $\mathbf{CP}^3$ branched along a singular  quartic surface that is bimeromorphic to an elliptic ruled surface.
We will see that the  branch elliptic curve of the minitwistor space is isomorphic to   the base elliptic curve of the branch quartic surface.
We also show that the isomorphism class of the branch elliptic curve uniquely determines the complex structure of the minitwistor space, and that the moduli space of our minitwistor spaces can be identified with the configuration space of different 4 points on $\mathbf{RP}^1\simeq S^1$ modulo the natural $\rm{PSL}(2,\mathbf R)$-action (Theorem \ref{thm-main2}).
In particular, our minitwistor spaces constitute a 1-dimensional moduli space, which contrasts to the case of the Gibbons-Hawking metrics and the LeBrun metrics where the moduli spaces are single point.

When investigating a twistor space, twistor lines  are of fundamental significance.
The images of  twistor lines into minitwistor space (by the quotient map) are  important as well and are called {\em minitwistor lines}.
In Section 3, we investigate minitwistor lines in our minitwistor spaces.
We prove that general minitwistor line is a nodal anticanonical curve of the minitwistor space, and that the natural morphism from a twistor line to the minitwistor line gives the normalization of the nodal curve.
Thus the situation is quite different from the case of the Gibbons-Hawking metrics and the LeBrun metrics, since in these two cases,   general minitwistor lines are smooth rational curves which are biholomorphic images of  twistor lines.
Geometrically, the appearance of the singularity of our minitwistor lines corresponds to the fact that for a general twistor line, there exists a unique $\mathbf C^*$-orbit intersecting  the twistor line {\em twice}. 
Finally we give an account why such situation occur in our (mini)twistor spaces  (Lemma \ref{lemma-rs2}), comparing that of LeBrun in \cite{LB91}.

The author would like to express his gratitude to Professors Shin Nayatani and Takashi Nitta for useful conversations. Also, he would like to thank Professor Akira Fujiki for his kind advise.

\vspace{1mm}
{\bf Notations.}
If $Z$ is a twistor space of a self-dual 4-manifold, $F$ denotes the canonical square root of the anticanonical line bundle $-K_Z$.
Tensor product of line bundles is denoted additively.

\section{The structure of minitwistor spaces and their moduli space}
First we recall the main results of \cite{H-bitan} which determine  global structure of the moduli space of self-dual metrics on $3\mathbf{CP}^2=\mathbf{CP}^2\#\mathbf{CP}^2\#\mathbf{CP}^2$ satisfying particular conditions.

\begin{prop}\label{prop-moduli1}
(\cite{H-bitan})
Let $g$ be a self-dual metric on $3\mathbf{CP}^2$ satisfying the following three properties:
(i) the scalar curvature of $g$ is  positive,
(ii) $g$ admits a non-zero killing field,
(iii) $g$ is not conformal to the self-dual metrics constructed by LeBrun \cite{LB91}.
Let $Z$ be the twistor space of $[g]$.
Then there is a commutative diagram of holomorphic maps
\begin{equation}\label{cd2}
\xymatrix{
   Z \ar@{->}[r]^{\mu}\ar@{->}[d]_{\Phi}  & Z_0  \ar@{->}[dl]^{\Phi_0}\\
   \mathbf{CP}^3 & \\
}
 \end{equation}
where $\Phi$ is a map associated to a linear system $|F|$ on  $Z$, $\Phi_0:Z_0\to\mathbf{CP}^3$ is a double covering whose branch locus is a quartic surface $B$ with ordinary nodes, $\mu$ is a small resolution of the corresponding ordinary nodes of $Z_0$.
Moreover, the defining equation of $B$ is given by
\begin{equation}\label{eqn-B}
\{y_2y_3+Q(y_0,y_1)\}^2-y_0y_1(y_0+y_1)(y_0-a y_1)=0,
\end{equation}
where $a$ is a positive real number and $Q(y_0,y_1)$ is a homogeneous quadratic polynomial with real coefficients satisfying the following condition:

\begin{quote}
{\em($\ast$)} as an equation on $\mathbf{CP}^1=\{(y_0,y_1)\}$, the quartic equation $$Q(y_0,y_1)^2-y_0y_1(y_0+y_1)(y_0-a y_1)=0$$ has a unique real double root.
\end{quote}

Conversely, for any quartic surface (\ref{eqn-B}) with $Q$ and $a>0$ satisfying the condition {\em($\ast$)}, the double covering of $\mathbf{CP}^3$ branched along $B$ admits a small resolution such that the resulting manifold is a twistor space of a self-dual metric on $3\mathbf{CP}^3$ satisfying  (i), (ii) and (iii).
\end{prop}

Apriori the real double root in the condition $(\ast)$ belongs to one of the two intervals $(-1,0)$ and $(a,\infty)$ on which the quartic $y_0y_1(y_0+y_1)(y_0-a y_1)$ is positive.
But as showed in \cite[\S 5.1]{H-bitan} we can always suppose that the double root belongs to the latter interval $(a,\infty)$ by applying a real projective transformation with respect to $(y_0,y_1)$, and in the following we always suppose this.
We also recall that the quartic surface $B$ defined by \eqref{eqn-B} satisfying the condition ($\ast$) has exactly three singular points
$$
P_{\infty}:=(0,0,0,1),\,\,\ol{P}_{\infty}:=(0,0,1,0),\,\,
P_0:=(\lambda_0,1,0,0),
$$
where $(\lambda_0,1)\in\mathbf{CP}^1$ is the real double root (so that $\lambda_0>a$ by the above normalization).

The killing field appeared in Proposition \ref{prop-moduli1} generates an isometric ${\rm{U}}(1)$-action of the self-dual metric.
This ${\rm{U}}(1)$-action naturally lifts and gives a holomorphic ${\rm{U}}(1)$-action on the twistor space $Z$.
Taking the complexification of the last ${\rm{U}}(1)$-action, we obtain a holomorphic $\mathbf C^*$-action on $Z$. 
This $\mathbf C^*$-action then descends on $\mathbf{CP}^3=\mathbf PH^0(F)^{\vee}$, which was shown to be of the form \cite[Prop. 2.1]{H-bitan}
\begin{equation}\label{eqn-action1}
(y_0,y_1,y_2,y_3)\longmapsto (y_0,y_1,ty_2,t^{-1}y_3),\,\,\,t\in\mathbf C^*.
\end{equation}
Of course this $\mathbf C^*$-action leaves the quartic surface $B$ invariant.
Note that any orbit of this action is contained in a plane belonging to the pencil $\langle y_0,y_1\rangle$, and that the closure of general orbits is a conic in these planes.
On the other hand, the anti-holomorphic involution of $\mathbf{CP}^3$ naturally induced from the real structure on $Z$ is explicitly given by, again as shown in \cite[Prop. 2.1]{H-bitan},
\begin{equation}\label{eqn-realstr}
(y_0,y_1,y_2,y_3)\longmapsto (\ol{y}_0,\ol{y}_1,\ol{y}_3,\ol{y}_2),
\end{equation}
where we are using the homogeneous coordinate in Proposition \ref{prop-moduli1}.
Of course, \eqref{eqn-action1} and \eqref{eqn-realstr} commute.

With these preliminary results, we begin to investigate  quotient spaces of the twistor spaces with respect to the $\mathbf C^*$-action.
First we consider a quotient space of the $\mathbf C^*$-action \eqref{eqn-action1} on $\mathbf{CP}^3$.
The rational map $\psi:\mathbf{CP}^3\to\mathbf{CP}^3$ defined by
\begin{equation}\label{quotient1}
\psi:(y_0,y_1,y_2,y_3)\longmapsto (z_0,z_1,z_2,z_3)=(y_0^2,y_1^2,y_0y_1,y_2y_3),
\end{equation}
which is the rational map associated to a linear system formed by $\mathbf C^*$-invariant quadratic polynomials, can be regarded as a quotient map of the $\mathbf C^*$-action, since general fibers of the map is the closure of general orbits.
The indeterminacy locus of $\psi$ consists of the two points $P_{\infty}$ and $\ol P_{\infty}$, which constitute a conjugate pair of points.
The image of $\psi$ is easily seen to be a quadratic cone $\ol{\Sigma}_2:=\{z_0z_1=z_2^2\}$
which has $(0,0,0,1)$ as the vertex.
Of course, $\ol{\Sigma}_2\backslash\{(0,0,0,1)\}$ is isomorphic to the total space of $\mathscr O(2)$, where the isomorphism is explicitly given by
\begin{align}\label{isom67}
\ol{\Sigma}_2\backslash\{(0,0,0,1)\}\ni
(z_0,z_1,z_2,z_3)\longmapsto (u,\zeta)=\left(\frac{z_2}{z_1},\frac{z_3}{z_1}\right)\in\mathscr O(2),
\end{align}
where $u$ is an affine coordinate on the base space of the bundle $\mathscr O(2)$ and $\zeta$ is a fiber coordinate valid on there.
Then by \eqref{eqn-B}, \eqref{eqn-action1}, \eqref{quotient1} and \eqref{isom67}, it is immediate to see the following

\begin{lemma}\label{lemma-imB}
Let $\psi$ be as in \eqref{quotient1} and $B$ the quartic surface in Proposition \ref{prop-moduli1}. Then under the isomorphism \eqref{isom67}, the image $\mathscr B:=\psi(B)$ is explicitly given by
\begin{equation}\label{cb}
\{\zeta+Q(u,1)\}^2-u(u+1)(u-a)=0.
\end{equation}
\end{lemma}

It is obvious from \eqref{cb} that the projection $\mathscr B\to\mathbf{CP}^1$ (given by $(u,\zeta)\mapsto u$) is a double covering
which has $u=0,-1, a$ as  simple branch points.
Further,
since \eqref{cb} is an equation taking values in the line bundle $\mathscr O(4)$, 
$u=\infty$ is also a simple branched point.
Also it is also obvious that these are all branch points.
Therefore, $\mathscr B$ is a smooth elliptic curve whose complex structure is determined by $a$.
Also we note that $\mathscr B$ does not go through the node of $\ol{\Sigma}_2$, and it belongs to the anticanonical class on $\ol{\Sigma}_2$.

The following result describes a structure of  quotient spaces of the twistor spaces in Proposition \ref{prop-moduli1} with respect to the $\mathbf C^*$-action:

\begin{thm}\label{thm-main1}
Let $g$, $Z$, $\Phi$ and $B$ as in Proposition \ref{prop-moduli1} and $\psi$, $\mathscr B$ as in Lemma \ref{lemma-imB}.
Then there exists a commutative diagram of meromorphic maps:
\begin{equation}\label{eqn-cd1}
\begin{CD}
Z@>{\Psi}>>\mathscr T\\
@V{\Phi}VV @VV{\phi}V\\
\mathbf{CP}^3@>{\psi}>>\ol{\Sigma}_2
\end{CD}
\end{equation}
where $\mathscr T$ is a normal rational surface, $\Psi$ is a surjective rational map, $\phi$ is a finite double covering map whose branch locus is the curve $\mathscr B$. 
Moreover, all fibers of $\Psi$ are $\mathbf C^*$-invariant and 
general fibers  are the closures of  orbits of the $\mathbf C^*$-actions on $Z$.
%
%
\end{thm}

By the last property, $\mathscr T$ can be regarded as a quotient space of the $\mathbf C^*$-action on $Z$.
Hence we call the normal rational surface $\mathscr T$ as a minitwistor space associated to the twistor space in Proposition \ref{prop-moduli1}.

\vspace{2mm}
\noindent Proof of Theorem \ref{thm-main1}.
Since $\psi$ has indeterminacy at $P_{\infty}$ and $\ol P_{\infty}$, 
the composition $\psi\circ\Phi$ also has indeterminacy along $\Phi^{-1}(P_{\infty})$ and $\Phi^{-1}(\ol P_{\infty})$.
As seen in \cite{H-bitan}, both of the last two sets are chains of 3 smooth rational curves in $Z$.
Let $Z'\to Z$ be  a sequence of blow-ups which resolves the indeterminacy of the rational map $\psi\circ\Phi$.
We may suppose that the image of the exceptional divisors are contained in $\Phi^{-1}(P_{\infty})\cup\Phi^{-1}(\ol P_{\infty})$.
Let $\Phi_0:Z_0\ra\mathbf{CP}^3$ be the double covering whose branch is $B$ as before.
If $x\in\ol{\Sigma}_2\backslash\mathscr B$, then $\psi^{-1}(x)$ is the closure of a $\mathbf C^*$-invariant conics which are not contained in $B$.
Hence $\Phi_0^{-1}(\psi^{-1}(x))$ splits into the closure of two orbits of the $\mathbf C^*$-action on $Z_0$.
On the other hand, if $x\in \mathscr B$, then $\psi^{-1}(x)$ is a $\mathbf C^*$-invariant conic contained in $B$.
Hence $\Phi_0^{-1}(\psi^{-1}(x))$ is biholomorphic to the conic $\psi^{-1}(x)$.
Therefore if $Z'\to\mathscr T\to\ol{\Sigma}_2$ is a Stein factorization of the morphism $Z'\to\ol{\Sigma}_2$, the latter $\mathscr T\to\ol{\Sigma}_2$ is a double covering whose branch is exactly $\mathscr B$.
Hence we obtained the commutative diagram \eqref{eqn-cd1}.
The statement about fibers of $\Psi$ is obvious from the above argument.
Further,  the singular locus of $\mathscr T$ is exactly the pre-image of the node of $\ol{\Sigma}_2$ since the branch curve $\mathscr B$ does not go through the node.
Hence $\mathscr T$ is normal.
Finally, the rationality of $\mathscr T$ is an immediate consequence of the fact that the pre-images of the lines on the cone $\ol{\Sigma}_2$ gives a pencil of rational curves on $\mathscr T$.
\proofend

\vsp
Since the $\mathbf C^*$-action on the twistor space is compatible with the real structure,
the minitwistor spaces also have real structures.
It is explicitly described as follows:

\begin{prop}\label{prop-rs}
Let $\mathscr T$ be the minitwistor space in Theorem \ref{thm-main1} and $\mathscr B\subset\ol{\Sigma}_2$ the branch elliptic curve of the double covering $\phi:\mathscr T\ra\ol{\Sigma}_2$.
Consider a real structure on $\ol{\Sigma}_2$ given by
\begin{equation}\label{eqn-rs55}
(u,\zeta)\longmapsto (\ol{u},\ol{\zeta}),
\end{equation}
on the complement $\mathscr O(2)$ of the node of  $\,\ol{\Sigma}_2$.
Then  $\mathscr B$ is invariant under this real structure, and the real structure on $\mathscr T$ covers the real structure \eqref{eqn-rs55} through the double covering $\mathscr T\ra\ol{\Sigma}_2$.
\end{prop} 

Proposition \ref{prop-rs} can be readily deduced 
by using \eqref{eqn-realstr}--\eqref{isom67}
 and we omit  the proof.

Another realization of the minitwistor space $\mathscr T$ is given by the following

\begin{prop}\label{cor-mt}
Let $\mathscr T$ be the minitwistor space in Theorem \ref{thm-main1}.
Then as a complex surface, $\mathscr T$ is obtained from $\mathbf{CP}^1\times\mathbf{CP}^1$ in the following way:
Let $p:\mathbf{CP}^1\times\mathbf{CP}^1 \ra\mathbf{CP}^1$ be (any) one of the projections, and $A_1$ and $A_2$ any two different sections of $p$ whose self-intersection numbers are zero.
Locate 4 points on $A_1\cup A_2$ in such a way that 2 points are on $A_1$ and the remaining  2 points are on $A_2$, and that the image of the 4 points under $p$ is equivalent to $\{-1,0,\infty,a\}$ under a projective transformation of $\mathbf{CP}^1$.
Next blowup $\mathbf{CP}^1\times\mathbf{CP}^1$ at these 4 points.
Then the strict transforms of $A_1$ and $A_2$ become $(-2)$-curves.
Finally blow down these two curves, to obtain a surface having two ordinary nodes.
This surface is biholomorphic to the minitwistor space $\mathscr T$.
\end{prop}

\noindent Proof.
Let $\nu:\Sigma_2\ra\ol{\Sigma}_2$ be the minimal resolution and write $\mathscr B'=\nu^{-1}(\mathscr B)$ which is isomorphic to $\mathscr B$.
Let $\mathscr T'\ra\Sigma_2$ be the double covering branched along $\mathscr B'$.
Then  $\mathscr T'$  is the minimal resolution of $\mathscr T$.
Consider the composition $\mathscr T'\ra\Sigma_2\ra\mathbf{CP}^1$, where $\Sigma_2\ra\mathbf{CP}^1$ is a projection of a ruling.
Since $\mathscr B'$ is 2 to 1 over $\mathbf{CP}^1$, general fiber of the above composition map is $\mathbf{CP}^1$.
Further, since $\mathscr B'$ has 4 branched points, the composition map has precisely 4 singular fibers, all of which are two $(-1)$ curves intersecting transversally.
If we choose four $(-1)$-curves among eight ones
in such a way that just two of them intersect one of the exceptional curves of the minimal resolution $\mathscr T'\ra\mathscr T$, and that the other two of them intersect another exceptional curve of $\mathscr T'\ra\mathscr T$,
 and if we  blow them down, then we obtain a (relatively) minimal surface which must be $\mathbf{CP}^1\times\mathbf{CP}^1$.
This implies the claim of the proposition.
\proofend

\vsp
We note that although Proposition \ref{cor-mt} gives an explicit construction of the minitwistor space as a complex surface, its real structure can never be obtained through this construction.
More precisely, 
the blowing-down $\mathscr T'\ra\mathbf{CP}^1\times\mathbf{CP}^1$ in the above proof does not preserve the real structure.
This can be seen, by going back to the twistor space, as follows.
Consider  singular fibers of $\mathscr T'\ra\mathbf{CP}^1$ in the above proof, which are pairs of $(-1)$-curves intersecting transversally at a point.
Then each of these singular fibers is the image of a reducible member of the linear system $|\Phi^*\mathscr O(1)|=|F|$, where $\Phi:Z\ra\mathbf{CP}^3$ is the generically 2 to 1 covering as in Proposition \ref{prop-moduli1}.
Namely,  the $(-1)$-curves are the images of the irreducible components, by the quotient map.
Since the real structure of $Z$ exchanges the irreducible components,
the two $(-1)$-curves in $\mathscr T'$ must be a conjugate pair.
Since the blowing-down  $\mathscr T'\ra \mathbf{CP}^1\times\mathbf{CP}^1$ contracts just one of the $(-1)$-curves for each reducible fiber, it cannot preserve the real structure.

\vsp
The complex structure of our twistor spaces in Proposition \ref{prop-moduli1} depends not only on $a>0$ but also on the coefficients of $Q(y_0,y_1)$ in the defining equation \eqref{eqn-B} of the branch quartic surface $B$.
We next show that the complex structure of the minitwistor spaces does not depend on $Q(y_0,y_1)$. Namely we show the following 

\begin{prop}\label{prop-moduli2}
Let $\mathscr T$ be the minitwistor space in Theorem \ref{thm-main1}.
Then the complex structure of $\mathscr T$ is uniquely determined by $a>0$ in the equation \eqref{eqn-B}. 
In other words, the complex structure of $\mathscr T$ does not depend on the quadratic polynomial $Q(y_0,y_1)$ in \eqref{eqn-B}. 
\end{prop}

\noindent
Proof. 
Fix $a>0$ and let $Q_1=Q_1(y_0,y_1)$ and $Q_2=Q_2(y_0,y_1)$ be two real homogeneous quadratic polynomials satisfying the condition ($\ast$) in Proposition \ref{prop-moduli1}.
Let $B_1$ and $B_2$ be the quartic surfaces determined by $(Q_1,a)$ and $(Q_2,a)$ by the equation \eqref{eqn-B} respectively.
Then we can write $Q_1(u,1)-Q_2(u,1)=d_0+d_1u+d_2u^2$ for $d_0,d_1,d_2\in\mathbf R$.
Using these $d_0,d_1,d_2$ we consider a map 
\begin{align}\label{isom23}
(u,\zeta)\mapsto (u,\zeta+d_0+d_1u+d_2u^2).
\end{align}
Viewing $(u,\zeta)$ as a holomorphic coordinate on the total space of $\mathscr O(2)$ as in the proof of Theorem \ref{thm-main1}, this map is easily seen to be a holomorphic automorphism of the Hirzebruch surface $\Sigma_2$.
Moreover, by \eqref{cb}, the automorphism \eqref{isom23} maps $\mathscr B_1$ to $\mathscr B_2$, where $\mathscr B_1$ and $\mathscr B_2$ are the images of $B_1$ and $B_2$ under the quotient map from $\mathbf{CP}^3$ to $\ol{\Sigma}_2$.
Thus we have concretely obtained an isomorphism of the pair $(\ol{\Sigma}_2,\mathscr B_1)$ and $(\ol{\Sigma}_2,\mathscr B_2)$.
Thus the double cover $\mathscr T_1$ and $\mathscr T_2$ whose branches are  $\mathscr B_1$ and $\mathscr B_2$ respectively are mutually biholomorphic, as desired.
\proofend

\vsp
Note that the isomorphism \eqref{isom23} between the pair $(\Sigma_2,\mathscr B_1)$ and $(\Sigma_2,\mathscr B_2)$ given in the above proof commutes with the real structure $(u,\zeta)\mapsto(\ol{u},\ol{\zeta})$, since $d_0, d_1$ and $d_2$ in \eqref{isom23}  are real.
Thus the minitwistor space $\mathscr T$ is uniquely determined by $a$ not only as a complex surface but also as a complex surface with real structure.

\vsp
By Proposition \ref{prop-moduli2} we can determine the moduli space of our minitwistor spaces. 
Let $\mathscr M$ be the moduli space of isomorphic classes of twistor spaces in Proposition \ref{prop-moduli1}.
As showed in \cite{H-bitan}, $\mathscr M$ is naturally identified with $\mathbf R^3/G$, where $G$ is a reflection of $\mathbf R^3$ having 2-dimensional fixed locus. 
Let $\mathscr N$ be the moduli space of isomorphic classes of  the associated minitwistor spaces, where the isomorphism is required to commute with the real structures.
We have a natural surjective map $\mathscr M\ra \mathscr N$ sending each isomorphic class of twistor space $Z$ to the isomorphic class of minitwistor space $\mathscr T$.
Then it is immediate from Proposition \ref{prop-moduli2} to obtain the following
\begin{thm}\label{thm-main2}
Let $\mathscr N$ be the moduli space of isomorphic classes of our minitwistor spaces as explained above.
Then $\mathscr N$ is naturally identified with the configuration space of different 4 points on a circle, divided by the usual {\rm{\,PSL}}$(2,\mathbf R)$-action on the circle.
\end{thm}
In particular  our minitwistor space has a non-trivial moduli, which contrasts to the known examples such as Gibbons-Hawking's \cite{GH78} and LeBrun's \cite{LB91,LB93}, where in these two cases the minitwistor spaces are the total space of $\mathscr O(2)$ and a quadratic surface $\mathbf{CP}^1\times\mathbf{CP}^1$  respectively and therefore
do not deform, although the corresponding self-dual (or hyperK\"ahler) metrics on 4-manifolds constitute non-trivial moduli spaces.

\section{Description of minitwistor lines}
As showed in the previous section, our minitwistor space $\mathscr T$ is a rational surface with two ordinary double points.
In this section we investigate minitwistor lines in $\mathscr T$; namely the images of twistor lines by the (rational) quotient map $\Psi:Z\ra\mathscr T$.
We investigate these minitwistor lines by using the diagram \eqref{eqn-cd1}.

A basic fact about twistor lines in our twistor space $Z$ was that, 
the image of general twistor line by the map $\Phi:Z\ra\mathbf{CP}^3$ is a very special kind of conic, called {\em touching conic}, meaning that the conic is tangent to the branch quartic surface $B$ at any intersection points which consist of 4 points in general
\cite[Def. 3.1 and Prop. 3.2]{H-bitan}.
Hence we first study the images of these touching conics by the (rational) quotient map $\psi:\mathbf{CP}^3\ra\ol{\Sigma}_2$:

\begin{lemma}\label{lemma-mtl1}
Let $\psi:\mathbf{CP}^3\ra\ol{\Sigma}_2$ be as in Lemma \ref{lemma-imB}.
Then the image  of general conics in $\mathbf{CP}^3$ under $\psi$ are  anticanonical curves on $\ol{\Sigma}_2$ with a unique node.
Further, this is true even for general touching conics of $B$, and their images  are  nodal anticanonical curves which touch the smooth anticanonical curve $\mathscr B$ at 4 points.
\end{lemma}

\noindent
Proof. Since any conic in $\mathbf{CP}^3$ is contained in some plane, we first study the restriction of $\psi$ onto a general plane. 
Since general orbits of our $\mathbf C^*$-action \eqref{eqn-action1} are conics, the restriction $\psi|_H$ is 2 to 1 for general plane $H$.
Further, by elementary calculations, we can readily see that $\psi|_H$ can be identified with a quotient map of $H=\mathbf{CP}^2$ by a reflection with respect to some line in $H$, where the line is exactly the set of tangents points of $\mathbf C^*$-orbits.
Further, the unique isolated fixed point of the reflection is mapped to the node of $\ol{\Sigma}_2$.
We say that a conic in a plane $H$ is symmetric if it is invariant under the reflection.
Then among the complete linear system $|\mathscr O(2)|$ on $H$, symmetric conics in $H$ form a codimension 2 linear subsystem.
It is easily seen that if a conic $C$ is not  symmetric, its image  is an anticanonical curve in $\ol{\Sigma}_2$ which has a unique node corresponding to the pair of intersection points of $C$ and its image in $H$ by the reflection which are not on the line.
(In contrasts, the image of symmetric conic becomes linearly equivalent to the branch curve of the map $H\to\ol{\Sigma}_2$.)
Thus we have seen that the image of a conic $C$ by $\psi$ is a nodal anticanonical curve in $\ol{\Sigma}_2$, as long as $C$ is not symmetric.
This shows the first claim of the proposition.

In order to show that the claim is still true for general touching conics of $B$, it suffices to show that for a plane smooth quartic $B_H$ on $H$ and any one of the 63 one-dimensional families of touching conics of $B_H$ (cf. \cite[Prop. 3.10]{H-bitan}), there exists no real line in $H$ for which all of the touching conics in the family become symmetric.
This can be proved as follows.
Let $\mathscr C$ be any one of the families of touching conics of $B_H$ and
suppose that there is a line $l$ in $H$ for which all  members of $\mathscr C$ are symmetric.
Then by \cite[Lemma 3.9 (ii)]{H-bitan} there are precisely 6 reducible members of $\mathscr C$, all of which are of course  pairs of two bitangents.
Since these reducible members must be symmetric with respect to $l$, the  intersection points of each pair of bitangents must  be on $l$.
Set $A=\Phi^{-1}_H(l)$, where $\Phi_H:S_H\ra H $ is the double cover of $H$ whose branch is $B_H$.
Then since $\Phi_H^*\mathscr O(1)\simeq -K_{S_H}$, $A$ is an anticanonical curve of $S_H$. We show that $A$ is irreducible.
For any irreducible member $C\in \mathscr C$, the inverse image $\Phi_H^{-1}(C)$ splits into a sum of two smooth rational curves $F_1$ and $F_2$ satisfying $F_1^2=F_2^2=0$ on $S_H$ (\cite[Lemma 3.8]{H-bitan}).
On $F_1$ and $F_2$, $\Phi_H$ is isomorphic to their images.
 Let $F$ be any one of $F_1$ and $F_2$ and consider the pencil $|F|$ generated by $F$.
Then members of $\mathscr C$ are exactly  the images of those of $|F|$ under the double covering map $\Phi_H$.
Because $c_1^2(S)=2$,
$|F|$ has precisely 6 reducible members and all of the members are mapped biholomorphically to a pair of bitangents that is a  reducible member of $\mathscr C$.
Therefore since $l$ goes through all the double points of the 6 reducible members of $\mathscr C$, $A=\Phi_H^{-1}(l)$ goes through 
every 6 points that are double points of the reducible members of the pencil $|F|$.
On the other hand, by considering the intersection of $A$ with any one of the reducible members of $|F|$, it immediately follows that  $F\cdot A=2$ on $S_H$, and that $A$ is not contained in any fiber of $h$.
Therefore, if $h:S_H\ra\mathbf{CP}^1$ denotes the morphism associated to the pencil $|F|$, the restriction $h|_A$ is 2 to 1  and  the 6 points  must be branch points of $h|_A$.
Hence $A$ must be irreducible, the 6 points on $S_H$ must be smooth points of $A$,  and  the 6 points  must be branch points of $h|_A$.
This means that the geometric genus of $A$ must be at least two.
This is a contradiction since $A$ is an anticanonical curve of the complex surface $S_H$.
Thus we have shown that there exists no line on $H$ with respect to which all members of $\mathscr C$ are symmetric.

Combined with what we have proved in the first paragraph of this proof, it follows that $\psi(C)$ is an anticanonical curve having a unique node, for a general touching conic $C$.
Since $\psi|_H$ is locally isomorphic outside the symmetric line on $H$, it follows that $\psi(C)$ still touches the image $\mathscr B=\psi(B)$ at four points.
Thus we have proved all the claims of the proposition.
\proofend

\vsp

Using Lemma \ref{lemma-mtl1} we show the following

\begin{prop}\label{prop-mtl2}
Let $\Psi:Z\ra\mathscr T$ be the (rational) quotient map by the $\mathbf C^*$-action on $Z$ as in Theorem \ref{thm-main1}.
Then the image of a general twistor line in $Z$ under $\Psi$ is a real anticanonical curve of $\mathscr T$ which has a unique node.
\end{prop}
In particular, general minitwistor lines in our minitwistor space $\mathscr T$ are not smooth.
This contrasts the case for LeBrun's metrics, since in LeBrun twistor spaces, since in LeBrun's case, minitwistor space is $\mathbf{CP}^1\times\mathbf{CP}^1$ and a general minitwistor line is a real (irreducible) curve of bidegree $(1,1)$, so that always non-singular.

\vsp
\noindent
Proof of Proposition \ref{prop-mtl2}.
As is already mentioned, $C=\Phi(L)$ is a touching conic of $B$ for a general twistor line $L$.
By Lemma \ref{lemma-mtl1} the image $\Gamma:=\psi(C)$ is a nodal anticanonical curve of $\ol{\Sigma}_2$.
Then by the diagram \eqref{eqn-cd1} the minitwistor line $\mathscr L:=\Psi(L)$ is an irreducible component of $\phi^{-1}(\Gamma)$.
As before let $\mu:\mathscr T'\ra\mathscr T$ and $\nu:\Sigma_2\ra\ol{\Sigma}_2$ be the minimal resolutions of  $\mathscr T$ and $\ol{\Sigma}_2$ respectively.
Let $\phi':\mathscr T'\to\Sigma_2$ be the natural lift of $\phi$.
Define  a line in $\mathbf{CP}^3$ by $l_{\infty}:=\{y_0=y_1=0\}$ in the coordinate of Proposition \ref{prop-moduli1}.
Then $l_{\infty}$ is exactly the fiber of $\psi:\mathbf{CP}^3\to\ol{\Sigma}_2$ over the node.
The branch locus of $\phi'$ is  a smooth anticanonical curve $\mathscr B'=\nu^{-1}(\mathscr B)$.
If $L$ is chosen so as to satisfy $\Phi(L)\cap l_{\infty}=\emptyset$, then $\Gamma=\Psi(\Phi(L))$ does not go through the node.
Hence $\Gamma':=\nu^{-1}(\Gamma)$ is  a nodal anticanonical curve of $\Sigma_2$ which is tangent to the branch curve $\mathscr B'$ at 4 points.
To prove the proposition, we have to look at irreducible components of 
the curve $(\phi')^{-1}(\Gamma')$.
It is immediate to see that $(\phi')^*(-K_{\Sigma_2})\simeq -2K_{\mathscr T'}$.
Moreover, since $\Gamma'$ is tangent to the branch curve $\mathscr B'$ at every intersection points, $(\phi')^{-1}(\Gamma')$ splits into two irreducible curves $\mathscr L_1$ and $\mathscr L_2$.
There are two possible situations:
(a) both $\mathscr L_1$ and $\mathscr L_2$ are smooth and the morphisms $\mathscr L_1\ra \Gamma'$ and $\mathscr L_2\ra \Gamma'$ (which are the restrictions of $\phi'$) are the normalizations of the nodal curve $\Gamma'$; or
(b) both $\mathscr L_1$ and $\mathscr L_2$ remain nodal curves and the morphisms $\mathscr L_1\ra \Gamma'$ and $\mathscr L_2\ra \Gamma'$ are isomorphic.
We now show that (a) cannot occur for general twistor lines by contradiction.
To this end, recall first that $\mathscr T'$ is realized as 4 points blown-up of  $\mathbf{CP}^1\times\mathbf{CP}^1$ as in Proposition \ref{cor-mt}.
In particular we have $(-K_{\mathscr T'})^2=4$ on $\mathscr T'$.
On the other hand, as is seen in the proof of Proposition \ref{cor-mt}, the restriction of the projection $\Sigma_2\ra\mathbf{CP}^1$ onto $\mathscr B'$ has 4 branch points and consequently the composition $\mathscr T'\ra\Sigma_2\ra\mathbf{CP}^1$ has precisely 4 singular fibers, all of which are the sum of two smooth rational curves intersecting transversally.
Let $E_i+E_i'$, $1\leq i\le 4$, be these 4 reducible fibers.
Then the blowing-down $\alpha:\mathscr T'\ra\mathbf{CP}^1\times\mathbf{CP}^1$ is obtained by  appropriately choosing one of $E_i$ and $E_i'$ for each $1\leq i\leq 4$ and then blowing them down.
After possible renaming, we suppose that $E_i$, $1\leq i\leq 4$, are blown-down by $\alpha$.
Then we  can write 
\begin{align}
\label{ls3}
\mathscr L_i=\alpha^*\mathscr O(a_i,b_i)-\sum_{j=1}^4n_{ij}E_j, \,\,i=1,2.
\end{align}
Since $\mathscr L_1+\mathscr L_2=-2K_{\mathscr T'}$, we have $a_1+a_2=b_1+b_2=4$ and $n_{1j}+n_{2j}=2$ for $1\leq j\leq 4$.
Moreover obviously we have $(a_i,b_i)=(1,3), (2,2)$ or $(3,1)$ for $i=1,2$.
Now we show that all $n_{ij}$ must be 1.
To see this, recall that the linear system $|F|$ on the twistor space has precisely 4 reducible members $\{D_i+\ol{D}_i\}_{i=1}^4$ and that all of them are $\mathbf C^*$-invariant.
Since  inverse images of fibers  of $\mathscr T\ra\mathbf{CP}^1$ by the quotient map $\Psi$ are $\mathbf C^*$-invariant members of $|F|$,
it follows that the 4 reducible fibers $E_i+E_i'$, $1\le i\le 4$, are the images of $D_i+\ol{D}_i$.
On the other hand, because $D_i\cdot L=\ol{D}_i\cdot L=1$,  general twistor lines intersect transversally with both of  $D_i$ and $\ol{D}_i$.
Hence $\mathscr L=\Psi(L)$ intersects both $E_i$ and $\ol{E}_i$ ($1\le i\le 4)$ for general $L$.
Thus combining with $n_{1j}+n_{2j}=2$, we have $n_{ij}=1$ for all $i$ and $j$.
Once this is proved, it readily follows $(a_1,b_1)=(a_2,b_2)=(2,2)$, since the blown-up 4 points of $\alpha$ are located in the way  described in Proposition \ref{cor-mt}, and in particular there are two sections of the projection with self-intersection zero on which 2 of the 4 blown-up points lie.
Thus we have 
\begin{align}\label{52}
\mathscr L_i=\alpha^*\mathscr O(2,2)-E_1-E_2-E_3-E_4 \,\,\text{ for }\,\,i=1,2.
\end{align}
Therefore the above item (a) cannot occur and (b) must hold for a general $L$.
Moreover, it is now obvious from \eqref{52} that $\mathscr L_1$ and $\mathscr L_2$ are anticanonical curves of $\mathscr T'$.
The reality of minitwistor lines is clear since the quotient map $\Psi$ preserve the real structure.
Thus we obtain all  claims of the proposition.
\proofend

\vsp
Finally we give another proof of the property that general minitwistor line in $\mathscr T$ has a node (Proposition \ref{prop-mtl2}), and compare the case of LeBrun twistor spaces

\begin{lemma}\label{lemma-rs2}
Consider the natural real structure on $\mathscr T$ which is induced from that on $Z$ (cf. Proposition \ref{prop-rs}).
Then the real locus on $\mathscr T$ consists of two disjoint 2-dimensional spheres.
Moreover, exactly one of the sphere parametrizes $\mathbf C^*$-orbits (in $Z$) whose closures are $\mathbf C^*$-invariant twistor lines.
\end{lemma}

\noindent
Proof. 
As in Proposition \ref{prop-rs}, the real structure on the total space of $\mathscr O(2)$ (which is the smooth locus of $\ol{\Sigma}_2$) is given by $(u,\zeta)\mapsto (\ol{u},\ol{\zeta})$.
Therefore the real locus of $\ol{\Sigma}_2$ consists of the closure of the set $\{(u,\zeta)\set u\in\mathbf R, \zeta\in\mathbf R\}$, which form a pinched torus, where the pinched point is the node of $\ol{\Sigma}_2$.
As is already seen, the branch locus $\mathscr B$ of the double covering $\phi:\mathscr T\ra\ol{\Sigma}_2$ is defined by
\begin{equation}\label{eqn-cb2}
\{\zeta+Q(u,1)\}^2-u(u+1)(u-a)=0.
\end{equation}
Hence the real locus on $\mathscr B$ is a union of the two sets given by
\begin{align}\label{rl1}
\left\{(u,\zeta)\in\mathscr O(2)\set -1\le u\le 0,\,\,\zeta=\pm\sqrt{u(u+1)(u-a)}\right\}
\end{align}
and 
\begin{align}\label{rl2}
\left\{(u,\zeta)\in\mathscr O(2)\set a\le u\le \infty,\,\,\zeta=\pm\sqrt{u(u+1)(u-a)}\right\}
\end{align}
where in the last condition we regard $\zeta=0$ if $u=\infty$.
The sets \eqref{rl1} and \eqref{rl2} are smooth circles in $\mathscr B$.
The sign of the left-hand side of \eqref{eqn-cb2} changes across these circles.
Since the double covering map $\mathscr T\to\ol{\Sigma}_2$ preserves the real structure, the real locus of $\mathscr T$ lies over the real locus of $\ol{\Sigma}_2$.
These mean that the real locus of $\mathscr T$ is either the inverse image of the two closed disks bounded by the circles \eqref{rl1} and \eqref{rl2}, or the inverse image of the complement of the last two disks.
But since the two points over the node of $\ol{\Sigma}_2$ (which is clearly outside the two circles) are a conjugate pair of points, the former must hold.
These two double covers of the closed disks are smooth spheres.

Next we see that exactly one of the two spheres parametrizes $\mathbf C^*$-orbits in $Z$ whose closures are $\mathbf C^*$-invariant twistor lines.
From the above description, the two spheres are over two intervals $[-1,0]$ and $[a,\infty]$ respectively, and every corresponding $\mathbf C^*$-orbits lie over a $\mathbf C^*$-invariant planes (determined by $u$).
Then as is shown in  \cite[Proposition 5.22]{H-bitan}, if $-1\le u\le 0$, then every real $\mathbf C^*$-orbits lying on the plane $y_0=uy_1$ must be an image of $\mathbf C^*$-invariant twistor lines.
Thus the sphere in $\mathscr T$ lying over $[-1,0]$ parametrizes $\mathbf C^*$-invariant twistor lines.
On the other hand, real orbits lying on a plane $y_0=uy_1$ with $u\ge a$ are not the image of twistor lines \cite[Proposition 5.22]{H-bitan}
This proves all the claims of the lemma.
\proofend

\vsp

By using Lemma \ref{lemma-rs2} we now give another explanation as to why general minitwistor lines in the minitwistor space $\mathscr T$ become singular.
Let $\mathscr T_2^{\sigma}$ and  $\mathscr T_4^{\sigma}$ be the connected components of the real locus $\mathscr T^{\sigma}$ of $\mathscr T$, where the former and latter lie over the interval $[-1,0]$ and $[a,\infty]$ in $\mathbf{RP}^1\subset\mathbf{CP}^1$ respectively.
(The subscripts $2$ and $4$ come from the notations in \cite{H-bitan}, where we wrote $I_2=[-1,0]$ and $I_4=[a,\infty]$.)
As above, both $\mathscr T_2^{\sigma}$ and  $\mathscr T_4^{\sigma}$ are 2-spheres smoothly embedded in $\mathscr T$.
As explained in the final part of the proof of Lemma \ref{lemma-rs2}, $\mathscr T_2^{\sigma}$ parametrizes $\mathbf C^*$-orbits in $Z$ whose closures are $\mathbf C^*$-invariant twistor lines,
while $\mathbf C^*$-orbits parametrized by $\mathscr T_4^{\sigma}$ are not (contained in) twistor lines.
Let $O\in\mathscr T^{\sigma}_4$ be any point and think $O$ as a real $\mathbf C^*$-orbit in $Z$.
Consider a twistor line $L\subset Z$ which intersects $O$.
Then by what we have explained above, $O$ is not contained in $L$ and 
it follows by reality that $L\cap O$ consists of even number of points.
The last number is 2 since $O$ is contained in some $S\in |F|$ and $S\cdot L=2$.
Then because the orbit map $\Psi:Z\ra \mathscr T$ identifies these conjugate pair of points, the image $\Psi(L)$ in $\mathscr T$ must have a singular point at $O\in\mathscr T$.
This is precisely the node of minitwistor line as stated in Proposition \ref{prop-mtl2}.
For each $O\in\mathscr T^{\sigma}_4$, there are obviously 2-dimensional family of twistor lines intersecting $L$.
Moreover, $\mathscr T^{\sigma}_4$ to which $O$ belongs, is also real 2-dimensional.
Thus there are real 4-dimensional family of twistor lines intersecting real orbit in $\mathscr T^{\sigma}_4$.
This means that the image of general twistor line must have a node.

\vsp
In contrast to the situation described in Lemma \ref{lemma-rs2}, the real locus of the minitwistor space of LeBrun twistor spaces on $n\mathbf{CP}^2$ consists of a unique sphere, and it parametrizes  $\mathbf C^*$-orbits whose closures are $\mathbf C^*$-invariant twistor lines.
This  is a reason why the image of general twistor lines by the orbit map is non-singular for LeBrun's twistor spaces.

\small
\vspace{13mm}
\hspace{5.5cm}
$\begin{array}{l}
\mbox{Department of Mathematics}\\
\mbox{Graduate School of Science and Engineering}\\
\mbox{Tokyo Institute of Technology}\\
\mbox{2-12-1, O-okayama, Meguro, 152-8551, JAPAN}\\
\mbox{{\tt {honda@math.titech.ac.jp}}}
\end{array}$

\end{document}